\documentclass[11pt]{article}
\usepackage{amssymb}
\usepackage{color}

\normalbaselines

\newtheorem{Theorem}{Theorem}[section]

\newtheorem{Definition}[Theorem]{Definition}

\newtheorem{Proposition}[Theorem]{Proposition}

\newtheorem{Lemma}[Theorem]{Lemma}

\newtheorem{Corollary}[Theorem]{Corollary}

\newtheorem{Remark}[Theorem]{Remark}

\newcommand{\RR}{{{\rm I} \kern -.15em {\rm R} }}

\newcommand{\C}{{{\rm l} \kern -.42em {\rm C} }}

\newcommand{\nat}{{{\rm I} \kern -.15em {\rm N} }}

\newcommand{\be}{\begin{equation}}
\newcommand{\ee}{\end{equation}}
\newcommand{\beq}{\begin{eqnarray}}
\newcommand{\eeq}{\end{eqnarray}}
\newcommand{\beqs}{\begin{eqnarray*}}
\newcommand{\eeqs}{\end{eqnarray*}}
\newcommand{\bt}{\begin{Theorem}}
\newcommand{\et}{\end{Theorem}}
\newcommand{\br}{\begin{Remark}}
\newcommand{\er}{\end{Remark}}
\newcommand{\bc}{\begin{Corollary}}
\newcommand{\ec}{\end{Corollary}}
\newcommand{\bl}{\begin{Lemma}}
\newcommand{\el}{\end{Lemma}}
\newcommand{\bd}{\begin{definition}}
\newcommand{\ed}{\end{definition}}

\title{Stability results for second--order evolution equations with\\ memory
and
switching time--delay}
\author{
{\sc Cristina Pignotti}
\\Dipartimento di Ingegneria e Scienze dell'Informazione e Matematica\\
 Universit\`{a} di L'Aquila\\
Via Vetoio, Loc. Coppito, 67010 L'Aquila Italy}

\date{}

\begin{document}

\textwidth=160 mm

\textheight=225mm

\parindent=8mm

\frenchspacing

\maketitle

\begin{abstract}
It is well-known
that wave--type equations with memory, under appropriate assumptions on the memory kernel, are uniformly exponentially stable. On the other hand, time delay effects may destroy this behavior. Here, we consider the stabilization problem for second-order evolution equations with memory and  intermittent delay feedback. We show that, under suitable assumptions involving the delay feedback coefficient and the memory kernel, asymptotic or exponential stability are still preserved. In particular, asymptotic stability is guaranteed if the delay feedback coefficient belongs to $L^1(0, +\infty )$
and the time intervals where the delay feedback is off are sufficiently large.
\end{abstract}

\vspace{5 mm}

\def\qed{\hbox{\hskip 6pt\vrule width6pt
height7pt
depth1pt  \hskip1pt}\bigskip}

%% {\bf 2000 Mathematics Subject Classification:}
%% 35L05, 93D15

 %%{\bf Keywords and Phrases:}  wave equation,  delay feedbacks, stabilization

\section{Introduction}
\label{pbform}\hspace{5mm}

\setcounter{equation}{0}

In this paper we study the stability properties of a viscoelastic model for second--order evolution equations.
In particular, we analyze a model combining memory damping and on--off time delay feedback, namely the time delay feedback is intermittently present.

It is well--known (see e.g. \cite{GRP, ACS}) that, under appropriate assumptions on the memory kernel, wave--type equations with viscoelastic damping are exponentially stable, i.e. the energy of all solutions is exponentially decaying to zero.
On the other hand, time  delay
effects
appear
in many applications and practical problems and it is by now well--known that even an arbitrarily small delay in the feedback may destabilize a system which is uniformly exponentially stable in absence of delay.
For some examples of this destabilizing effect of time delays we refer to \cite{Datko, DLP, NPSicon06, XYL}.

We want to show that under suitable conditions involving  the delay feedback coefficient  and the memory kernel, the system is asymptotically stable or exponentially stable, in spite of the presence of the time delay term.

Stability results for second--order evolution equations with intermittent damping
are first studied by Haraux, Martinez and Vancostenoble \cite{HMV}, without any time delay term. They consider a problem with intermittent on--off or with positive--negative damping,
and show that, under appropriate conditions,  the {\em good} behavior of the system in the time intervals where  only the standard dissipation, i.e. the damping with the right sign, is present prevails  over the {\em bad} behavior where
the damping is no present or it is present with the wrong sign, i.e. as anti--damping.
Thus, asymptotic/exponential stability results are obtained. See also \cite{Genni} for a nonlinear extension.

Recently  Nicaise and Pignotti \cite{ADE2012, JDDE14} considered
second--order evolution equations with intermittent delay feedback.
More precisely, in the studied models, when the (destabilizing) delay term is no present,
a not--delayed damping acts. Under suitable assumptions, stability results are obtained.
Other results for wave equations with intermittent delay feedback have been obtained, in 1-dimension, in \cite{Gugat}, \cite{GT} and \cite{ANP2012} by using a completely different approach based on  the D'Alembert formula
and only  for particular choices of the time delay.

In the present paper, the good behavior in the time intervals where the delay feedback is no present is ensured by a viscoelastic damping.

Let $H$ be a real Hilbert space and let $A:{\mathcal D}(A)\rightarrow H$
be a positive self--adjoint  operator with a compact inverse in $H.$ Denote by $V:={\mathcal D}(A^{\frac 1 2})$ the domain of
$A^{\frac 1 2}.$

Let us consider the problem
 \begin{eqnarray}
& &u_{tt}(x,t) +A u (x,t)-\int_0^\infty \mu(s) A u(x, t-s) ds +b(t) u_t(x, t-\tau) =0\quad t>0,\label{1.1}\\
& & u(x,t) = 0 \quad \mbox{\rm on}\quad \partial\Omega\times (0,+\infty ),\label{CB}\\
& &u(x,t)=u_0(x,t)\quad \mbox{\rm in}\quad \Omega\times (-\infty, 0];\label{1.2}
\end{eqnarray}

where  the initial datum $u_0$ belongs to a suitable space, the constant $\tau >0$ is the time delay, and
the  memory kernel $\mu :[0,+\infty)\rightarrow [0,+\infty)$

satisfies

i) $\mu\in C^1 (\RR^+ )\cap L^1(\RR^+ );$

ii) $\mu (0)=\mu_0>0;$

iii) $\int_0^{+\infty} \mu (t) dt=\tilde \mu <1;$

iv) $\mu^{\prime} (t)\le -\delta \mu (t), \quad \mbox{for some}\ \ \delta >0.$

Moreover, the function $b(\cdot )\in L^\infty_{loc}(0,+\infty )$ is a function which is zero intermittently. That is,
we assume that for all $n\in\nat$ there exists $t_n>0$ with $t_n<t_{n+1}$
and such that
$$
\begin{array}{l}
b(t)=0\quad  \forall\ t\in I_{2n}=[t_{2n},t_{2n+1}),\\
\vert b(t)\vert < b_{2n+1}\ne 0
\quad \forall \  t\in I_{2n+1}=[t_{2n+1},t_{2n+2}).
\end{array}
$$

 Also, denoting by $T_n$ the length of the interval $I_n,$ that is
\begin{equation}\label{Tn}
T_n=t_{n+1}-t_n,\quad n\in \nat\,,
\end{equation}
we assume that 
\begin{equation}\label{T2n}
\tau\leq T_{2n}, \quad  \forall\  n\in\nat\,.
\end{equation}

\noindent
We know that the above problem is exponentially stable for $ b(t)\equiv 0 $ (see e.g.
\cite{GRP}).

Without memory (i.e. $\mu\equiv 0$) and undelayed damping, namely $\tau=0,$ exponential stability results have been obtained,  when $b(t)\equiv b >0$  (see \cite{Rauch} and also, for various generalizations, \cite{BLR, Komornikbook, Lag83, LT, Lions, liu, zuazua}). However, when a time delay is present, the standard frictional damping generates unstabilities. Then, a stabilizing term, in this case the memory damping, has to be included in order to gain stability.

Models with both viscoelastic damping and time delay feedback have been studied by several recent papers.
The stability properties of the wave equation with memory and time delay have been first studied by Kirane and Said-Houari \cite{KSH}, in the case of finite memory. However, in their model, an extra standard dissipative damping is added in order to contrast the destabilizing effect of the time delay term.

On the contrary, more recent papers show that the viscoelastic damping is sufficient to contrast the delay term and then to ensure exponential stability under a suitable {\em smallness} condition on the time delay coefficient.
The stabilization problem for wave--type equations  with infinite memory and time delay feedback  has been studied
by Guesmia in \cite{Guesmia} by using a suitable Lyapunov functional and by Alabau--Boussouira, Nicaise and Pignotti \cite{AlNP2015} by using a perturbation approach for delay problems first introduced in \cite{pignotti}.
We refer also to Day and Yang
 \cite{DY} for the same kind problem  in the case of finite memory.
 In these papers the authors prove exponential stability results if the (constant) coefficient of the delay damping is sufficiently small. These stability results can be easily extended to
a variable coefficient $b(\cdot )\in L^\infty (0, +\infty )$ under a suitable {\em smallness} assumption on the  $L^\infty-$ norm of $b(\cdot )\,.$

In the present paper, instead,
we  show that an asymptotic stability result holds without smallness conditions related to the $L^\infty-$norm.
On the other hand, we are restricted to consider on--off delay term in order to
obtain our results. We show, in particular, that asymptotic stability is ensured if $b(\cdot )\in L^1(0, +\infty )$ and the length of the time intervals where the delay feedback is off, i.e. $b\equiv 0,$ is sufficiently large. We also give stability results, under analogous assumptions, for a problem with on--off anti--damping instead of a time delay feedback.

The paper is organized as follows. In section \ref{well}    a well--posedness result for the abstract system is proved. In section \ref{st} we obtain,
asymptotic and exponential stability results for the abstract model
 under suitable conditions. We also give the stability results for the problem with on--off anti--damping.
Finally, in section  \ref{esempi}, we illustrate our abstract results by some concrete applications.

\section{Well-posedness \label{well}}

\hspace{5mm}

\setcounter{equation}{0}

In this section we will give a well-posedness results for problem
(\ref{1.1})--(\ref{1.2}).

For the existence result let us consider the problem

 \begin{eqnarray}
& &u_{tt}(x,t) +A u (x,t)-\int_0^\infty \mu(s) A u(x, t-s) ds  = f(t)\quad t>0,\label{NH.1}\\
& & u(x,t) = 0 \quad \mbox{\rm on}\quad \partial\Omega\times (0,+\infty ),\label{NH.B}\\
& &u(x,t)=u_0(x,t)\quad \mbox{\rm in}\quad \Omega\times (-\infty, 0].\label{NH.2}
\end{eqnarray}

As in Dafermos \cite{Dafermos}, let us introduce the new variable
\begin{equation}\label{eta}
\eta^t(x,s):=u(x,t)-u(x,t-s).
\end{equation}

Then, problem (\ref{NH.1})--(\ref{NH.2}) may be rewritten as

\begin{eqnarray}
& &u_{tt}(x,t)= -(1-\tilde \mu)A u (x,t)-
\int_0^\infty \mu (s) A \eta^t(x,s) ds\nonumber\\
& &\hspace{5 cm}
 +f(t)\quad \mbox{\rm in}\ \Omega\times
(0,+\infty)\label{e1d}\\
& & \eta_t^t(x,s)=-\eta^t_s(x,s)+u_t(x,t)\quad \mbox{\rm in}\ \Omega\times
(0,+\infty)\times (0,+\infty ),\label{e2d}\\
& &u (x,t) =0\quad \mbox{\rm on}\ \partial\Omega\times
(0,+\infty)\label{e3d}\\
& &\eta^t (x,s) =0\quad \mbox{\rm in}\ \partial\Omega\times
(0,+\infty), \ t\ge 0,\label{e4d}\\
& &u(x,0)=u_0(x)\quad \mbox{\rm and}\quad u_t(x,0)=u_1(x)\quad \hbox{\rm
in}\ \Omega,\label{e5d}\\
& & \eta^0(x,s)=\eta_0(x,s) \quad \mbox{\rm in}\ \Omega\times
(0,+\infty), \label{e6d}
\end{eqnarray}
where
\begin{equation}\label{datiinizd}
\begin{array}{l}
u_0(x)=u_0(x,0), \quad x\in\Omega,\\
u_1(x)=\frac {\partial u_0}{\partial t}(x,t)\vert_{t=0},\quad x\in\Omega,\\
\eta_0(x,s)=u_0(x,0)-u_0(x,-s),\quad x\in\Omega,\  s\in (0,+\infty).
\end{array}
\end{equation}
Set
$L^2_{\mu}((0, \infty); V)$  the Hilbert space
of $V-$ valued functions on $(0,+\infty),$
endowed with the inner product
$$\langle \varphi, \psi\rangle_{L^2_{\mu}((0, \infty);V)}=
\int_0^\infty \mu (s)\Vert A^{1/2} \varphi (s)\Vert^2 ds\,.
$$

\noindent
Denote by ${\mathcal H}$
the Hilbert space

$${\mathcal H}=
V\times H\times L^2_{\mu}((0, \infty); V),$$
equipped
  with the inner product

\begin{equation}\label{innerd}
\begin{array}{l}
\left\langle
\left (
\begin{array}{l}
u\\
v\\
w
\end{array}
\right ),\left (
\begin{array}{l}
\tilde u\\
\tilde v\\
\tilde w
\end{array}
\right )
\right\rangle_{\mathcal H}
:= \displaystyle{
 (1-\tilde\mu )\langle A^{1/2} u, A^{1/2}\tilde u \rangle_H + \langle v, \tilde v \rangle_H +
 \int_0^\infty \mu (s)\langle A^{1/2} w, A^{1/2}\tilde w \rangle_H ds }
\end{array}
\end{equation}

Let us recall the following well--posedness result (see \cite{GRP}).

\begin{Definition}\label{Pata}
Set $I=[0,T],$ for $T>0,$ and let $f\in L^1(I, H).$ A function
$U:= (u, u_t, \eta )\in {\mathcal H}$ is a solution of problem $(\ref{e1d})-(\ref{datiinizd})$
in the interval $I,$ with initial data $U(0)=U_0= (u_0, u_1, \eta_0)\in {\mathcal H},$
provided
$$\langle u_{tt}, \tilde v\rangle_H=-(1-\tilde\mu )\langle A^{1/2} u, A^{1/2}\tilde v
\rangle_H-\int_0^\infty \mu (s)
\langle A^{1/2}\eta (s), A^{1/2}\tilde v\rangle_H ds  +\langle f,\tilde v\rangle_H\,;$$
$$\int_0^\infty \mu (s) \langle  \eta_t(s)+\eta_s(s),A\tilde\eta (s)  \rangle_H ds
=\int_0^\infty \mu(s) \langle u_t, A \tilde\eta (s)\rangle_H ds\,;$$
for all $\tilde v\in V$ and $\tilde\eta \in L^2_{\mu}(\RR^+, {\mathcal D}(A)),$ and a.e. $t\in I.$
\end{Definition}

\begin{Theorem}\label{noto}
For given $T>0,$ problem $(\ref{e1d})-(\ref{datiinizd})$ has a unique solution $U$ in the time interval $I=[0,T],$ with initial datum $U_0.$
\end{Theorem}

Under these assumptions, we obtain the following result

\begin{Theorem}\label{texistence}
Under the above assumptions, for any $U_0\in {\mathcal H},$ the
system $(\ref{1.1})-(\ref{1.2})$ has a unique solution
 $U\in C([0,\infty); {\mathcal H})$.
\end{Theorem}
\noindent {\bf Proof.} We prove the existence and uniqueness result on the interval $[0,t_2];$ then the global result follows by translation (cfr. \cite{ADE2012}). First, in the interval $[0, t_1],$ since $b(t)= 0$ $\forall\ t\in [0, t_1),$ we can apply Theorem \ref{noto} with $f\equiv 0\,.$ Then we obtain a solution $U,$ in the sense of Definition \ref{Pata}, on the interval $[0, t_1].$ The situation is more delicate in the time interval $[t_1, t_2]$ where the delay feedback is present. In this case, we decompose the interval $[t_1, t_2]$ into the successive intervals $[t_1+j\tau, t_1+(j+1)\tau ),$ for $j=0, \dots, N,$ where $N$ is such that $t_1+(N+1)\tau\ge t_2\,.$
The last interval is then $[t_1+N\tau, t_2]\,.$
Now, look at the problem on the interval $[t_1, t_1+\tau]\,.$ Here $u_t(t-\tau)$ can be considered as a known function. Indeed, for $t\in [t_1, t_1+\tau]\,,$ then $t-\tau \in [0, t_1],$ and we know the solution $U$ on $[0, t_1]$ by the first step. Thus, problem
$(\ref{1.1})-(\ref{1.2})$ may be rewritten on $[t_1, t_1+\tau ]$ as

 \begin{eqnarray}
& &u_{tt}(x,t) +A u (x,t)-\int_0^\infty \mu(s) A u(x, t-s) ds =f_1(t)\quad t\in [t_1, t_1+\tau ],\label{W1.1}\\
& & u(x,t) = 0 \quad \mbox{\rm on}\quad \partial\Omega\times [t_1,t_1+\tau] ,\label{WCB}\\
& &u(x,t)=u_0^1(x,t)\quad \mbox{\rm in}\quad \Omega\times (-\infty, t_1];\label{W1.2}
\end{eqnarray}

where $f_1(t)=b(t)u_t(t-\tau)\,$ the initial datum is  $u_0^1(x,t)= u_0(x,t)$ in $\Omega\times (-\infty, 0]$ and $u_0^1(x,t)= u(x,t)$ in $\Omega\times [0, t_1]\,.$

Then we can apply once more Theorem \ref{noto} obtaining a unique solution $U$ on $[0, t_1+\tau )\,.$ Proceedings analogously in the successive time intervals $[t_1+j\tau, t_1+(j+1)\tau ),$ we obtain a solution on $[0, t_2]\,.$\qed

\section{Stability results
\label{st}}

\hspace{5mm}

\setcounter{equation}{0}

As in \cite{NPSicon06} we introduce the variable
\begin{equation}\label{zeta}
z(x,\rho ,t):=
u_t(x,t-\tau\rho),\quad x\in\Omega,\ \rho\in (0,1),\ t>0.
\end{equation}

\noindent
Using (\ref{eta}) and (\ref{zeta}) we can rewrite (\ref{1.1})--(\ref{1.2})
as

\begin{eqnarray}
& &u_{tt}(x,t)= -(1-\tilde \mu)A u (x,t)-
\int_0^\infty \mu (s)A\eta^t(x,s) ds\nonumber\\
& &\hspace{5 cm}
 -b(t)
z(x,1,t)\quad \mbox{\rm in}\ \Omega\times
(0,+\infty)\label{e1d2}\\
& & \eta_t^t(x,s)=-\eta^t_s(x,s)+u_t(x,t)\quad \mbox{\rm in}\ \Omega\times
(0,+\infty)\times (0,+\infty ),\label{e2d2}\\
& &\tau z_t(x,\rho ,t)+z_{\rho}(x,\rho ,t)=0\quad \mbox{\rm in}\ \Omega\times
(0,1)\times (0,+\infty ),\label{e2dzeta2}\\
& &u (x,t) =0\quad \mbox{\rm on}\ \partial\Omega\times
(0,+\infty)\label{e3d2}\\
& &\eta^t (x,s) =0\quad \mbox{\rm in}\ \partial\Omega\times
(0,+\infty), \ t\ge 0,\label{e4d2}\\
& &z(x,0,t)=u_t(x,t)\quad \mbox{\rm in}\ \Omega\times
(0,+\infty), \label{e4dzeta2}\\
& &u(x,0)=u_0(x)\quad \mbox{\rm and}\quad u_t(x,0)=u_1(x)\quad \hbox{\rm
in}\ \Omega,\label{e5d2}\\
& & \eta^0(x,s)=\eta_0(x,s) \quad \mbox{\rm in}\ \Omega\times
(0,+\infty), \label{e6d2}\\
& &z(x,\rho , 0)=z^0(x,-\tau\rho) \quad x\in\Omega, \ \rho\in (0,1),\label{e7d2}
\end{eqnarray}
where
\begin{equation}\label{datiinizd2}
\begin{array}{l}
u_0(x)=u_0(x,0), \quad x\in\Omega,\\
u_1(x)=\frac {\partial u_0}{\partial t}(x,t)\vert_{t=0},\quad x\in\Omega,\\
\eta_0(x,s)=u_0(x,0)-u_0(x,-s),\quad x\in\Omega,\  s\in (0,+\infty),\\
z^0(x,s)=\frac {\partial u_0}{\partial t}(x,s),\quad x\in\Omega,\ s\in (-\tau ,0).
\end{array}
\end{equation}

Let us now introduce the energy functional
\begin{equation}\label{energyd}
\begin{array}{l}
\displaystyle{
E(t)=E(u,t):=\frac 1 2 \Vert u_t(t)\Vert^2_H
+\frac {1-\tilde\mu} 2 \Vert u(t)\Vert^2_V}\\
\hspace{1 cm}\displaystyle{
+\frac 1 2 \int_0^{+\infty } \mu(s)\Vert A^{1/2}\eta^t(s)\Vert^2_H ds +\frac 12\int_{t-\tau}^t \vert b(s+\tau )\vert\Vert  u_t(s)\Vert^2_H ds \,.}
\end{array}
\end{equation}
Then,
$$E(t)= E_S(t)+ \frac 12\int_{t-\tau}^t \vert b(s+\tau )\vert\Vert  u_t(s)\Vert^2_H ds \,,$$
where $E_S(\cdot )$ denotes the standard energy for the wave equation with viscoelastic damping, i.e.
\begin{equation}\label{standardenergy}
E_S(t)=\frac 1 2 \Vert u_t(t)\Vert^2_H
+\frac {1-\tilde\mu} 2 \Vert u(t)\Vert^2_V
+\frac 1 2 \int_0^{+\infty } \mu(s)\Vert A^{1/2}\eta^t(s)\Vert^2_H ds\,.
\end{equation}

Let us now recall the following result proved in \cite{GRP} for wave equation (see e.g. \cite{Pata} for  the abstract case).

\begin{Theorem}\label{Pata}
Assume $b\equiv 0\,.$ Then, for every solution of problem $(\ref{1.1})-(\ref{1.2}),$ the energy $E_S(\cdot )$ is not increasing and
\begin{equation}\label{decay}
E_S^\prime (t)\le \frac 1 2 \int_0^\infty \mu^\prime (s)\Vert A^{1/2}\eta^t(s)\Vert^2_H ds\,.
\end{equation}
Moreover,
there are two positive constant $C,\alpha  ,$  $C>1, \alpha >0,$ depending only on $\Omega$ and on the memory kernel $\mu (\cdot ),$ such that for every solution of problem $(\ref{1.1})-(\ref{1.2})$ it results
\begin{equation}\label{exp}
E_S(t)\le Ce^{-\alpha t} E_S(0)\,.
\end{equation}
\end{Theorem}

Now, let $T_0$ be the time such that

\begin{equation}\label{timeweneed}
T_0:=\frac 1 {\alpha } \ln C\,,
\end{equation}
that is the time for which $Ce^{-\alpha T}=1\,.$

As an immediate application we have the following result.

\begin{Proposition}\label{obs}
Assume $T_{2n}>T_0\,.$ Then, there exists a constant $c_n\in (0,1)$ such that
\begin{equation}\label{stimaobs}
E_S(t_{2n+1})\le c_n E_S (t_{2n})\,,
\end{equation}
for any solution of problem $(\ref{1.1})-(\ref{1.2}).$
\end{Proposition}

\noindent {\bf Proof.}
Observe that in the time interval $I_{2n}=[t_{2n}, t_{2n+1}]$ the delay feedback is not present since $b\equiv 0\,.$
Then, from estimate $(\ref{exp})$ and the semigroup property, we deduce
$$E_S(t)\le C e^{-\alpha (t-t_{2n})}E_S(t_{2n}),\quad \forall\ t\in I_{2n}\,.$$
Then,
$$ E_S(t_{2n+1})\le C e^{-\alpha T_{2n}}E_S(t_{2n})\,.$$
If $T_{2n}$ is greater than $T_0,$ where $T_0$ is as in (\ref{timeweneed}), it is immediate to see that
\begin{equation}\label{cn}
c_n:=Ce^{-\alpha T_{2n}}\in (0,1)\,;
\end{equation}
then the claim is proved. \qed

\begin{Proposition}\label{intervallicattivi}
Assume $T_{2n}\ge\tau,\ \forall\ n\in \nat .$ Then,
\begin{equation}\label{stimacattiva}
E^\prime (t) \le  b_{2n+1} \Vert u_t(t)\Vert_H^2\,,\quad t\in I_{2n+1}= [t_{2n+1}, t_{2n+2}], \ \forall\ n\in \nat\,.
\end{equation}
for any solution of problem $(\ref{1.1})-(\ref{1.2}).$
\end{Proposition}

\noindent {\bf Proof.}
By differentiating the energy $E(\cdot ),$ we have
$$\begin{array}{l}
\displaystyle{
E^\prime (t)=\langle u_t(t), u_{tt}(t)\rangle_H+ (1-\tilde\mu ) \langle A^{1/2}u(t),
A^{1/2}u_t(t)\rangle_H}\\
\hspace{1 cm}\displaystyle{
+\int_0^\infty \mu (s)\langle A^{1/2}\eta^t(s), A^{1/2}\eta^t_t(s)\rangle_H ds
+\frac 1 2 \vert b(t+\tau )\vert \Vert u_t(t)\Vert_H^2-\frac 1 2 \vert b(t)\vert \Vert u_t(t-\tau )\Vert_H^2\,.}
\end{array}
$$

Then, by Green's formula and using the boundary condition,

\begin{equation}\label{cri10}
 \begin{array}{l}
\displaystyle{
E^\prime (t)=\langle u_t(t), u_{tt}(t)-(1-\tilde\mu )A u(t)\rangle_{V-V^\prime}}\\
\hspace{1 cm}\displaystyle{
+\int_0^\infty \mu (s)\langle A^{1/2}\eta^t(s), A^{1/2}u_t(t) -A^{1/2}\eta^t_s(s)\rangle_H ds}\\
\hspace{1 cm}\displaystyle{
+\frac 1 2 \vert b(t+\tau )\vert \Vert u_t(t)\Vert_H^2-\frac 1 2 \vert b(t)\vert \Vert u_t(t-\tau )\Vert_H^2\,.}
\end{array}
\end{equation}

By (\ref{cri10}), by using the equation (\ref{1.1}) and integrating by parts, we obtain
 $$
 \begin{array}{l}
\displaystyle{
E^\prime (t)=\langle u_t(t), \int_0^\infty \mu(s) A u(t-s) ds -\tilde\mu A u(t)-b(t)u_t(t-\tau) \rangle_{V-V^\prime}}\\
\hspace{1 cm}\displaystyle{
+\int_0^\infty \mu (s)\langle A^{1/2}\eta^t(s), A^{1/2}u_t(t)\rangle_H ds
+\frac 12\int_0^\infty\mu^\prime (s)
 \Vert A^{1/2}\eta^t(s)\Vert^2_H ds}\\
\hspace{1 cm}\displaystyle{
+\frac 1 2 \vert b(t+\tau )\vert \Vert u_t(t)\Vert_H^2-\frac 1 2 \vert b(t)\vert \Vert u_t(t-\tau )\Vert_H^2\,.}
\end{array}
$$
Finally, from Cauchy--Schwarz inequality,
$$E^\prime (s)\le \frac 1 2 \vert b(t)\vert \Vert u_t(t)\Vert_H^2+\frac 1 2 \vert b(t+\tau )\vert \Vert u_t(t)\Vert_H^2\le b_{2n+1}\Vert u_t(t)\Vert_H^2\,,$$
where in the last inequality we have used the fact that, since $T_{2n}\ge\tau,$ for every $n\in\nat\,,$ if $t\in I_{2n+1}$ then
$t+\tau\in I_{2n+1}\cap I_{2n+2}\,.$
Thus (\ref{stimacattiva}) is proved.\qed

\begin{Theorem}\label{CP1}
Assume $\mbox{\rm i), ii)}$ and $T_{2n}\ge \tau$ for all $n\in\nat.$ Moreover
assume $T_{2n}> T_0,$ for all $n\in\nat,$ where $T_0$ is the time defined in $(\ref{timeweneed}).$ Then, if
\begin{equation}\label{general}
\sum_{n=0}^\infty
\ln \left [
e^{2b_{2n+1}T_{2n+1}}
 (c_n+T_{2n+1} b_{2n+1} )\right ] = -\infty\,,
\end{equation}
the
system $(\ref{1.1})-(\ref{1.2})$ is asymptotically stable, that is any solution  $u$ of $(\ref{1.1})-(\ref{1.2})$ satisfies
$E_S(u,t)\rightarrow 0$ for $t\rightarrow +\infty\,.$
\end{Theorem}

\noindent {\bf Proof.}
Note that (\ref{stimacattiva}) implies
$$E^{\prime}(t)\le 2b_{2n+1} E(t),\quad t\in I_{2n+1}=[t_{2n+1},t_{2n+2}),\ n\in\nat.$$
Then we have
\begin{equation}\label{cri0}
E(t_{2n+2})\le e^{2b_{2n+1}T_{2n+1}}E(t_{2n+1}),
\quad \forall \ n\in\nat.
\end{equation}
From the definition of the energy $E$,
\begin{equation}\label{cri1}
E(t_{2n+1})=E_S(t_{2n+1})+\frac 1 2 \int_{t_{2n+1}-\tau }^{t_{2n+1}}\vert b(s+\tau )\vert\Vert  u_t(s)\Vert^2_H  ds\,.
\end{equation}
Note that, for $t\in [t_{2n+1}-\tau , t_{2n+1}),$ then $t+\tau \in [t_{2n+1}, t_{2n+1}+\tau )\subset I_{2n+1}\cup I_{2n+2}\,.$
Now, if $t+\tau \in I_{2n+2},$  then $b(t+\tau )=0.$ Otherwise, if $t+\tau \in I_{2n+1},$  then $\vert b(t+\tau)\vert\le b_{2n+1}\,.$
Then, from (\ref{cri1}) we deduce
\begin{equation}\label{cri2}
E(t_{2n+1})=E_S(t_{2n+1})+\frac 1 2  b_{2n+1} \int_{t_{2n+1}-\tau }^{\min ({t_{2n+2}-\tau }, t_{2n+1})}\Vert u_t(s)\Vert^2_H ds\,,
\end{equation}
since if  $t_{2n+1}> t_{2n+2}-\tau ,$ then  $b(t)=0$ for all $t\in [t_{2n+2}, t_{2n+1}+\tau )\subset [t_{2n+2}, t_{2n+3}).$

Then, since the energy $E_S(\cdot )$ is decreasing in the intervals $I_{2n},$
\begin{equation}\label{cri3}
E(t_{2n+1})\le E_S(t_{2n+1})+T_{2n+1} b_{2n+1} E_S(t_{2n+1}-\tau )
\le E_S(t_{2n+1})+T_{2n+1} b_{2n+1} E_S(t_{2n})
\,.
\end{equation}

Using (\ref{cri3}) in (\ref{cri0}), we deduce

\begin{equation}\label{cri4}
E_S(t_{2n+2})\le E(t_{2n+2})\le e^{2b_{2n+1}T_{2n+1}}
 (c_n+ T_{2n+1} b_{2n+1}
)E_S(t_{2n}),
\quad \forall \ n\in\nat,
\end{equation}
where we have used also the observability estimate (\ref{stimaobs}).
Iterating this procedure we arrive at

\begin{equation}\label{cri5}
E_S(t_{2n+2})\le\displaystyle{\Pi_{p=0}^n }e^{2b_{2p+1}T_{2p+1}}
 (c_p+T_{2p+1} b_{2p+1}
)E_S(0),
\quad \forall \ n\in\nat\,.
\end{equation}

Now observe that the standard energy $E_S(\cdot )$ is not decreasing in $(0,+\infty).$ However, it is decreasing
for $t\in [t_{2n}, t_{2n+1})$, when only the standard dissipative damping acts and so
\begin{equation}\label{Nice1}
E_S(t)\le E_S(t_{2n}),\quad \forall \; t\in [t_{2n}, t_{2n+1}).
\end{equation}
Moreover, from (\ref{cri3}), for $t\in [t_{2n+1}, t_{2n+2})$ it results
\begin{equation}\label{Nice2}
E_S(t)\le E(t)\le e^{2b_{2n+1}T_{2n+1}}(c_n+b_{2n+1}T_{2n+1}) E(t_{2n}),
\end{equation}
where in the second inequality we have used  (\ref{stimaobs}).

Then, we have asymptotic stability if

$$\displaystyle{
\Pi_{p=0}^n} e^{2b_{2p+1}T_{2p+1}}
 (c_p+T_{2p+1} b_{2p+1}) \longrightarrow 0,\quad \mbox{\rm for}\ n\rightarrow\infty,$$
 or equivalently
 $$\ln \Big [
\Pi_{p=0}^n e^{2b_{2p+1}T_{2p+1}}
 (c_p+T_{2p+1} b_{2p+1} )\Big ] \longrightarrow -\infty,\quad \mbox{\rm for}\ n\rightarrow\infty\,.$$
 This concludes the proof. \qed

\begin{Remark}\label{particular}
{\rm
In particular
 (\ref{general}) is verified if the following conditions are satisfied:
\begin{equation}\label{star1bis}
\sum_{n=0}^\infty b_{2n+1}T_{2n+1}<+\infty\quad \mbox{\rm and}\quad \sum_{n=0}^\infty \ln c_n=-\infty\,.
\end{equation}
Indeed, it is easy to see that
\begin{equation}\label{star1}
\sum_{n=0}^\infty b_{2n+1}T_{2n+1}<+\infty\quad \mbox{\rm and}\quad \sum_{n=0}^\infty \ln (c_n+b_{2n+1}T_{2n+1})=-\infty\,,
\end{equation}
imply (\ref{general}).
Now, observe also that $(\ref{star1})$ is equivalent to $(\ref{star1bis}).$
Indeed, if $(\ref{star1})$ holds true then also $(\ref{star1bis})$ is verified, since
$$\ln c_n <\ln (c_n+b_{2n+1}T_{2n+1})\,, \quad \forall\ n\in\nat\,.$$
Now assume that $(\ref{star1bis})$ holds true. From the first condition we have
\begin{equation}\label{conditnec}
b_{2n+1}T_{2n+1}\rightarrow 0\quad \mbox{\rm  for}\quad  n\rightarrow +\infty\,.
\end{equation}
Suppose by contradiction that $(\ref{star1})$ does not hold. This implies that
$$\sum_{n=0}^\infty \ln (c_n+b_{2n+1}T_{2n+1})>-\infty$$
and then, being the terms of the series definitely negative,
$$\sum_{n=0}^\infty \ln (c_n+b_{2n+1}T_{2n+1})\in (-\infty, 0)\,.$$
Therefore, it has to be
$\ln (c_n+b_{2n+1}T_{2n+1})\rightarrow 0$ or equivalently
$c_n+b_{2n+1}T_{2n+1}\rightarrow 1,$ as $n\rightarrow\infty \,.$
We conclude, by $(\ref{conditnec}),$ that
$$c_n+b_{2n+1}T_{2n+1} \sim c_n\,.$$
 Then  $(\ref{star1})$ has to be verified.

In particular, from $(\ref{star1bis}),$  we have stability if $b\in L^1(0,+\infty )$ and,
for instance,
the length of the {\em good} intervals $I_{2n}$ is greater than a fixed time $\bar T,$ $\bar T>T_0$ and $\bar T\ge\tau ,$
namely
$$T_{2n}\ge \bar T,\quad\forall\ n\in \nat\,.$$ Indeed, in this case there exists $\bar c\in (0,1)$ such that $0<c_n<\bar c\,.$ }
\end{Remark}

We now show that under additional assumptions on the coefficients $T_n, b_{2n+1}, c_n$ an exponential stability  result holds.

\begin{Theorem}\label{exp}
Assume $\mbox{\rm i), ii)}$
and that
\begin{equation}\label{E1}
T_{2n}=T^*\quad\forall\ n\in\nat ,
\end{equation}
with $T^*\ge \tau$ and $T^*>T_0,$ where the time $T_0$ is as in $(\ref{timeweneed}).$
Assume also that
\begin{equation}\label{E2}
T_{2n+1}=\tilde T\quad\forall\ n\in\nat .
\end{equation}
Moreover,  assume that
\begin{equation}\label{ASS1A}
\sup_{n\in\nat} \ e^{2 b_{2n+1}\tilde T}(c+b_{2n+1}T_{2n+1})=d<1,
\end{equation}
where $c_n=c,\ n\in\nat,$ is as in $(\ref{stimaobs})$. Then, there exist two positive constants $\gamma,\beta$
such that
\begin{equation}\label{expestimateA}
E_S(t)\le \gamma e^{-\beta t} E_S(0),\ \ t>0,
\end{equation}
for any solution of problem  $(\ref{1.1})-(\ref{1.2}).$
\end{Theorem}

\noindent {\bf Proof.} From (\ref{ASS1A}) and (\ref{cri5})
we obtain
$$E_S(T^*+\tilde T)\le d E(0),$$
and also
$$E_S(n(T^*+\tilde T))\le d^n E(0),\quad \forall n\in\nat.$$
Then, the energy satisfies an exponential estimate like
(\ref{expestimateA}) (see Lemma 1 of \cite{Gugat}).\qed

\begin{Remark}{\rm
In the assumptions of Theorem \ref{exp},
from (\ref{cri5}) we can see that exponential stability also holds if instead of
(\ref{ASS1A}) we assume
$$\exists
 n\in\nat\quad \mbox{\rm such}\ \mbox{\rm that}\quad
\displaystyle{\Pi_{p=k(n+1)}^{k(n+1)+n}} e^{2 b_{2p+1}\tilde T}(c+b_{2p+1}T_{2p+1})\le d<1,\ \ \forall\ k=0,1,2,\dots
$$
}\end{Remark}

\begin{Remark}\label{finitememory}{\rm
Analogous results could be obtained for the case of finite memory, namely if system (\ref{1.1})-(\ref{1.2}) is replaced by
 \begin{eqnarray}
& &u_{tt}(x,t) +A u (x,t)-\int_0^t \mu(s) A u(x, t-s) ds +b(t) u_t(x, t-\tau) =0\quad t>0,\label{F.1}\\
& & u(x,t) = 0 \quad \mbox{\rm on}\quad \partial\Omega\times (0,+\infty ),\label{FCB}\\
& &u(x,t)=u_0(x,t),\quad u_t(x,t)=u_1(x)(\mbox{\rm in}\quad \Omega;\label{F1.2}
\end{eqnarray}
with memory kernel $\mu (\cdot )$ and delay coefficient $b(t)$ satisfying the same assumptions that before.
Indeed, also for such a problem it is well-known that an exponential decay estimate holds on the time intervals where the delay feedback is null (see e.g. \cite{ACS}). Therefore, in such intervals an observability type estimate like (\ref{stimaobs}) is available if the length of the intervals is sufficiently large.
}\end{Remark}

\subsection{Stability under the restriction $T_{2n+1}\le\tau$\label{st1}}

\hspace{5mm}

Now, we assume that the length of the delay intervals is lower
than the time
delay, that is
\begin{equation}\label{rest}
T_{2n+1}\le\tau,\quad \forall n\in\nat\,.
\end{equation}

We look at the standard energy $E_S(\cdot).$
We can give the following estimates
on the time intervals $I_{2n}, I_{2n+1},$ $n\in\nat.$

\begin{Proposition}\label{CP2}
Assume $\mbox{\rm i),\ ii)}.$ Moreover assume $T_{2n+1}\le\tau$ and $T_{2n}\ge\tau\,.$
Then, for
 $t\in I_{2n+1},$
\begin{equation}\label{Pi1}
E_S^{\prime}(t)\le
b_{2n+1}E_S(t)+
b_{2n+1}E_S(t_{2n})\,.
\end{equation}
\end{Proposition}

\noindent{\bf Proof:} By differentiating $E_S(t)$ we get
\begin{eqnarray*}
E_S'(t)=b(t)\langle u_t(t), u_t(t-\tau )\rangle_H\,.
\end{eqnarray*}

Hence,
from ${\rm ii}),$
$$
E_S'(t)\le \frac {b_{2n+1}}2\Vert u_t(t)\Vert^2_H+\frac {b_{2n+1}}2\Vert u_t(t-\tau )\Vert^2_H\le b_{2n+1}E_S(t)+b_{2n+1}E_S(t-\tau )\,.
$$
Now, to conclude it suffices to observe that since $T_{2n+1}\le\tau$ and $T_{2n}\ge\tau\,,$ then for $t\in I_{2n+1}$ it is $t-\tau \in I_{2n}.$ince
Then, since $E_S(\cdot )$ is decreasing in $I_{2n}\,,$ estimate (\ref{Pi1}) is proved.
\qed

\begin{Theorem}\label{CP3}
Assume $\mbox{\rm i),\ ii)}.$ Moreover assume $T_{2n+1}\le\tau$ and $T_{2n}\ge\tau\,,$
$\forall\ n\in \nat\,.$
Then, if
\begin{equation}\label{star3}
\sum_{n=0}^\infty
\ln \left [
e^{b_{2n+1}T_{2n+1}}
 (c_n+1-e^{-b_{2n+1} T_{2n+1}}
  )\right ] = -\infty\,,
\end{equation}
the system $(\ref{1.1})-(\ref{1.2})$ is asymptotically stable, that is any solution  $u$ of $(\ref{1.1})-(\ref{1.2})$ satisfies
$E_S(u,t)\rightarrow 0$ for $t\rightarrow +\infty\,.$
\end{Theorem}

\noindent {\bf Proof.} For $t\in I_{2n+1}=[t_{2n+1}, t_{2n+2}),$
from estimate (\ref{Pi1}) we deduce
$$
E_S^{\prime}(t)\le e^{b_{2n+1}(t-t_{2n+1})}\Big [ E_S(t_{2n+1})+ \int_{t_{2n+1}}^t b_{2n+1}E_S(t_{2n}) e^{-b_{2n+1}(s-t_{2n+1})} ds\Big ]\,.
$$

Then we have
$$
E_S(t)\le e^{b_{2n+1}T_{2n+1} }E_S(t_{2n+1})+
e^{b_{2n+1}(t-t_{2n+1})}E_S(t_{2n})\Big ( 1-e^{-b_{2n+1}(t-t_{2n+1})}\Big )\,,
$$
and so
$$E_S(t)\le e^{b_{2n+1}T_{2n+1} }E_S(t_{2n+1})+E_S(t_{2n})e^{b_{2n+1}T_{2n+1} }-E_S(t_{2n})\,,$$
for $t\in I_{2n+1}=[t_{2n+1},t_{2n+2}),\ n\in\nat.$

In particular, by using estimate (\ref{stimaobs}), we obtain

$$E_S(t_{2n+2})\le e^{b_{2n+1} T_{2n+1}}[  c_n +1-e^{-b_{2n+1} T_{2n+1}}
]
E_S(t_{2n}),\quad n\in\nat\,,$$
and therefore
\begin{equation}\label{Pi2}
E_S(t_{2n+2})\le \Big (
\displaystyle{
\Pi_{p=0}^n}e^{b_{2p+1} T_{2p+1}} [ c_p +1 -e^{-b_{2p+1} T_{2p+1}}]
\Big )E_S(0)\,\end{equation}
Then, by (\ref{Pi2}), asymptotic stability occurs if
$$
 \displaystyle{
\Pi_{p=0}^n}e^{b_{2p+1} T_{2p+1}} [ c_p +1 -e^{-b_{2p+1} T_{2p+1}}]
 \rightarrow 0,\quad \mbox{\rm for}\ n\rightarrow\infty\,,
$$

or, equivalently, if

$$
 \sum_{n=0}^\infty
\ln \Big ( e^{b_{2n+1} T_{2n+1}} [ c_n +1 -e^{-b_{2n+1} T_{2n+1}}]
\Big ) \rightarrow -\infty,\quad \mbox{\rm for}\ n\rightarrow\infty\,.
$$
This concludes.
\qed

\begin{Remark}{\rm
Note that, in case of {\em bad} intervals $I_{2n+1}$ with length lower or equal than the time delay $\tau\,,$ the assumption (\ref{star3}) is a bit less restrictive than (\ref{general})\,.
Indeed, since $b_{2n+1}T_{2n+1}>0,$ it results
$$e^{b_{2n+1}T_{2n+1}}(c_n+1-e^{-b_{2n+1}T_{2n+1}})
<e^{2b_{2n+1}T_{2n+1}}(c_n+b_{2n+1}T_{2n+1}),\quad \forall n\in \nat\,.$$
For instance if $b_{2n+1}T_{2n+1}=1/4$ and $c_n= e^{-1/2}-1/4$ for every $n\in\nat,$
then
$$ e^{2b_{2n+1}T_{2n+1}}(c_n+b_{2n+1}T_{2n+1})=1, \quad \forall\ n\in \nat\ ,$$
and
$$e^{b_{2n+1}T_{2n+1}}(c_n+1-e^{-b_{2n+1}T_{2n+1}})=\alpha \in (0,1),\quad \forall\ n\in \nat\ .$$
Therefore (\ref{general}) does not hold while (\ref{star3}) is clearly satisfied.
}
\end{Remark}

\begin{Remark}{\rm Arguing as in Remark \ref{particular} we can show that
condition (\ref{star3}) is verified, in particular, if (\ref{star1bis}) holds true.}
\end{Remark}

Also in this case, under additional assumptions on the coefficients $T_n, b_{2n+1}, c_n$ an exponential stability  result holds.

\begin{Theorem}\label{exp2}
Assume $\mbox{\rm i), ii)}$
and that
$$
T_{2n}=T^*\quad\forall\ n\in\nat ,
$$
with $T^*\ge \tau$ and $T^*>T_0,$ where the time $T_0$ is as in $(\ref{timeweneed}).$
Assume also that
\begin{equation}\label{E2bis}
T_{2n+1}=\tilde T, \quad \mbox{\rm with}\ \ \tilde T\le\tau \quad\forall\ n\in\nat .
\end{equation}
Moreover,  assume that
\begin{equation}\label{ASS1Abis}
\sup_{n\in\nat} \ e^{ b_{2n+1}\tilde T}(c+1 -e^{-b_{2n+1}T_{2n+1}})=d<1,
\end{equation}
where $c_n=c,\ n\in\nat,$ is as in $(\ref{stimaobs})$. Then, there exist two positive constants $\gamma,\beta$
such that
\begin{equation}\label{expestimateAbis}
E_S(t)\le \gamma e^{-\beta t} E_S(0),\ \ t>0,
\end{equation}
for any solution of problem  $(\ref{1.1})-(\ref{1.2}).$
\end{Theorem}

\noindent {\bf Proof.} The proof is analogous to the one of Theorem \ref{exp}\,.\qed

\subsection{Problem with anti--damping\label{AD}}

\hspace{5mm}
With the same technics we can also deal with an intermittent anti--damping term.
More precisely, let us consider the problem
 \begin{eqnarray}
& &u_{tt}(x,t) +A u (x,t)-\int_0^\infty \mu(s) A u(x, t-s) ds -k(t) u_t(x, t) =0\quad t>0,\label{1.1A}\\
& & u(x,t) = 0 \quad \mbox{\rm on}\quad \partial\Omega\times (0,+\infty ),\label{CBA}\\
& &u(x,t)=u_0(x,t)\quad \mbox{\rm in}\quad \Omega\times (-\infty, 0];\label{1.2A}
\end{eqnarray}

where  the initial datum $u_0$ belongs to a suitable space, 
the  memory kernel $\mu :[0,+\infty)\rightarrow [0,+\infty)$ is as before
and
the function $k(\cdot )\in L^\infty_{loc}(0,+\infty )$ is a function which is zero intermittently. That is,
we assume that for all $n\in\nat$ there exists $t_n>0$ with $t_n<t_{n+1}$
and such that
$$
\begin{array}{l}
k(t)=0\quad  \forall\ t\in I_{2n}=[t_{2n},t_{2n+1}),\\
\vert k(t)\vert < k_{2n+1}\ne 0
\quad \forall \  t\in I_{2n+1}=[t_{2n+1},t_{2n+2}).
\end{array}
$$

In particular, $k$ may be positive, i.e. an anti--damping term (cfr. \cite{HMV}).
 As before, denote by $T_n$ the length of the interval $I_n,$ that is
$$
T_n=t_{n+1}-t_n,\quad n\in \nat\,.
$$
Now, the time delay is no present. Therefore, we deal with the standard energy (\ref{standardenergy}).

Of course Proposition \ref{obs} holds, which gives an observability estimate on the intervals $I_{2n}$ where the anti-damping is off.
On the time intervals $I_{2n+1}$ one can obtain the following estimate.

\begin{Proposition}\label{badAD}
For every solution of problem $(\ref{1.1A})-(\ref{1.2A}),$ 
$$E_S^\prime (t)\le  k_{2n+1}  \Vert u_t(t)\Vert_H^2, \quad t\in I_{2n+1}=[t_{2n+1}, t_{2n+2}], \ \ \forall \ n\in \nat\,.$$
\end{Proposition}

\begin{Theorem}\label{CP1AD}
Assume $\mbox{\rm i), ii)}.$ Moreover
assume $T_{2n}> T_0,$ for all $n\in\nat,$ where $T_0$ is the time defined in $(\ref{timeweneed}).$ Then, if
\begin{equation}\label{generalAD}
\sum_{n=0}^\infty
\ln \left (
e^{2k_{2n+1}T_{2n+1}}
 c_n\right ) = -\infty\,,
\end{equation}
the
system $(\ref{1.1A})-(\ref{1.2A})$ is asymptotically stable, that is any solution  $u$ of $(\ref{1.1A})-(\ref{1.2A})$ satisfies
$E_S(u,t)\rightarrow 0$ for $t\rightarrow +\infty\,.$
\end{Theorem}

\noindent {\bf Proof.}
Note that from Proposition \ref{badAD} we have 
$$E_S^{\prime}(t)\le 2b_{2n+1} E_S(t),\quad t\in I_{2n+1}=[t_{2n+1},t_{2n+2}),\ n\in\nat.$$
Then we have
\begin{equation}\label{cri0A}
E_S(t_{2n+2})\le e^{2k_{2n+1}T_{2n+1}}E_S(t_{2n+1}),
\quad \forall \ n\in\nat.
\end{equation}
Then, from estimate (\ref{stimaobs}),

\begin{equation}\label{cri4A}
E_S(t_{2n+2})\le e^{2k_{2n+1}T_{2n+1}}
 c_n
E_S(t_{2n}),
\quad \forall \ n\in\nat\,.
\end{equation}
Iterating this procedure we arrive at

\begin{equation}\label{cri5}
E_S(t_{2n+2})\le\displaystyle{\Pi_{p=0}^n }e^{2k_{2p+1}T_{2p+1}}
 c_pE_S(0),
\quad \forall \ n\in\nat\,.
\end{equation}
Then, we have asymptotic stability if

$$\displaystyle{
\Pi_{p=0}^n} e^{2k_{2p+1}T_{2p+1}}
 c_p \longrightarrow 0,\quad \mbox{\rm for}\ n\rightarrow\infty,$$
 or equivalently
 $$\ln \Big (
\Pi_{p=0}^n e^{2k_{2p+1}T_{2p+1}}
 c_p\Big ) \longrightarrow -\infty,\quad \mbox{\rm for}\ n\rightarrow\infty\,.$$
 This concludes the proof. \qed

\begin{Remark}\label{particularAD}
{\rm
In particular
 (\ref{generalAD}) is verified if the following conditions are satisfied:
\begin{equation}\label{star1bisAD}
\sum_{n=0}^\infty k_{2n+1}T_{2n+1}<+\infty\quad \mbox{\rm and}\quad \sum_{n=0}^\infty \ln c_n=-\infty\,.
\end{equation}
}
\end{Remark}

Under additional assumptions on the coefficients $T_n, k_{2n+1}, c_n$ an exponential stability  result holds.

\begin{Theorem}\label{expAD}
Assume $\mbox{\rm i), ii)}$
and that
\begin{equation}\label{E1AD}
T_{2n}=T^*\quad\forall\ n\in\nat ,
\end{equation}
with $T^*>T_0,$ where the time $T_0$ is as in $(\ref{timeweneed}).$
Assume also that
\begin{equation}\label{E2AD}
T_{2n+1}=\tilde T\quad\forall\ n\in\nat .
\end{equation}
Moreover,  assume that
\begin{equation}\label{ASS1AAD}
\sup_{n\in\nat} \ e^{2k_{2n+1}\tilde T}c=d<1,
\end{equation}
where $c_n=c,\ n\in\nat,$ is as in $(\ref{stimaobs})$. Then, there exist two positive constants $\gamma,\beta$
such that
\begin{equation}\label{expestimateAAD}
E_S(t)\le \gamma e^{-\beta t} E_S(0),\ \ t>0,
\end{equation}
for any solution of problem  $(\ref{1.1A})-(\ref{1.2A}).$
\end{Theorem}

\section{Examples\label
{esempi}}
\hspace{5mm}

\setcounter{equation}{0}

As a concrete example we can consider the wave equation with memory.
More precisely, let $\Omega\subset\RR^n$ be an open bounded domain   with a smooth  boundary
$\partial\Omega.$
Let us consider the initial boundary value  problem

\begin{eqnarray}
& &u_{tt}(x,t) -\Delta u (x,t)+
\int_0^\infty \mu (s)\Delta u(x,t-s) ds\nonumber \\
& &\hspace{5 cm}
 +
b(t) u_t(x,t-\tau )=0\quad \mbox{\rm in}\ \Omega\times
(0,+\infty)\label{1.1d}\\
& &u (x,t) =0\quad \mbox{\rm on}\ \partial\Omega\times
(0,+\infty)\label{1.2d}\\
& &u(x,t)=u_0(x, t)\quad \hbox{\rm
in}\ \Omega\times (-\infty, 0] \label{1.3d}
%%\\
%%u_t(x,t)=f_0(x,t) \quad \mbox{\rm in} \ \Omega\times (-\tau, 0),
%%\label{1.4d}
\end{eqnarray}

This problem enters into our previous framework, if we take
$H=L^2(\Omega)$ and the operator  $A$ defined by
$$A:{\mathcal D}(A)\rightarrow H\,:\,  u\rightarrow -\Delta u,$$
where ${\mathcal D}(A)=H^1_0(\Omega)\cap
H^2(\Omega).$

The operator $A$ is a self--adjoint and positive operator with a compact inverse in $H$
and is such that $V={\mathcal D}(A^{1/2})=H^1_0(\Omega).$

The energy functional is, in this case,

\begin{equation}\label{energyd}
\begin{array}{l}
\displaystyle{
E(t)=\frac 1 2 \int_{\Omega} u_t^2(x,t) dx
+\frac {1-\tilde\mu} 2 \int_{\Omega}\vert \nabla u(x,t)\vert^2 dx}\\
\hspace{1 cm}\displaystyle{
+\frac 1 2 \int_0^{+\infty } \int_{\Omega}\mu(s)\vert \nabla\eta^t(s)\vert^2 ds dx+\frac 1 2\int_{t-\tau}^t \int_{\Omega}\vert b(s+\tau )\vert u_t^2(x,s) ds dx.}
\end{array}
\end{equation}

Under the same conditions that before on the memory kernel $\mu (\cdot)$ and on the function $b(\cdot ),$ previous asymptotic/exponential stability results are valid.
The case $b$ constant has been studied in \cite{AlNP2015} by adapting a perturbative approach introduced in \cite{pignotti}.
In particular, we have proved that the exponential stability is preserved, in presence of the delay feedback, if the coefficient of this one is sufficiently small.
The choice $b$ constant was made only for the sake of clearness. The result in  \cite{AlNP2015} remains true if instead of $b$ constant we consider $b=b(t),$
under a suitable smallness condition on the $L^\infty-$norm of   $b(\cdot )\,.$
On the contrary here we give stability results without restrictions on the  $L^\infty-$norm of   $b(\cdot )\,,$
even if only for on--off $b(\cdot )\,.$

Our results also apply to Petrovsky  system with viscoelastic damping  with Dirichlet and Neumann boundary conditions:

\begin{eqnarray}
& &u_{tt}(x,t) +\Delta^2 u (x,t)-
\int_0^\infty \mu (s)\Delta^2 u(x,t-s) ds\nonumber \\
& &\hspace{5 cm}
 +
b(t) u_t(x,t-\tau )=0\quad \mbox{\rm in}\ \Omega\times
(0,+\infty)\label{1.1P}\\
& &u (x,t) =\frac {\partial u}{\partial\nu }=0\quad \mbox{\rm on}\ \partial\Omega\times
(0,+\infty)\label{1.2P}\\
& &u(x,t)=u_0(x, t)\quad \hbox{\rm
in}\ \Omega\times (-\infty, 0] \label{1.3P}
%%\\
%%u_t(x,t)=f_0(x,t) \quad \mbox{\rm in} \ \Omega\times (-\tau, 0),
%%\label{1.4d}
\end{eqnarray}

This problem enters into the previous abstract framework, if we take
$H=L^2(\Omega)$ and the operator  $A$ defined by
$$A:{\mathcal D}(A)\rightarrow H\,:\,  u\rightarrow \Delta^2 u,$$
where ${\mathcal D}(A)=H^2_0(\Omega)\cap
H^4(\Omega),$ with
$$H^2_0(\Omega )=\Big\{
v\in H^2(\Omega )\ :\ v=\frac {\partial v}{\partial\nu }=0 \ \ \mbox{\rm on}\ \ \partial\Omega\
\Big\}\,.$$
The operator $A$ is a self--adjoint and positive operator with a compact inverse in $H$
and is such that $V={\mathcal D}(A^{1/2})=H^2_0(\Omega).$

In this case, the energy functional becomes

\begin{equation}\label{energyP}
\begin{array}{l}
\displaystyle{
E(t)=\frac 1 2 \int_{\Omega} u_t^2(x,t) dx
+\frac {1-\tilde\mu} 2 \int_{\Omega}\vert \Delta u(x,t)\vert^2 dx}\\
\hspace{1 cm}\displaystyle{
+\frac 1 2 \int_0^{+\infty } \int_{\Omega}\mu(s)\vert \Delta\eta^t(s)\vert^2 ds dx+\frac 1 2\int_{t-\tau}^t \int_{\Omega}\vert b(s+\tau )\vert u_t^2(x,s) ds dx.}
\end{array}
\end{equation}

 Therefore, under the same conditions that before on the memory kernel $\mu (\cdot)$ and on the function $b(\cdot ),$ previous asymptotic/exponential stability results are valid.

\bigskip

 {\em E-mail address,}
\\
\quad Cristina Pignotti: \quad{\tt \bf
pignotti@univaq.it}

\end{document}